\documentstyle{amsppt}
\topmatter
\title{REARRANGEMENTS OF TRIGONOMETRIC SERIES AND TRIGONOMETRIC POLYNOMIALS}
\endtitle
\author{S.~V.~Konyagin}
\endauthor
\abstract
The paper is related to the following question of P.~L.~Ul'yanov:
is it true that for any $2\pi$-periodic continuous function $f$ there is 
a uniformly convergent rearrangement of its trigonometric Fourier series? 
In particular, we give an affirmative answer if the absolute values 
of Fourier coefficients of $f$ decrease. Also, we study a problem how 
to choose $m$ terms of a trigonometric polynomial of degree $n$
to make the uniform norm of their sum as small as possible.
\endabstract
\endtopmatter

\centerline{\S1. Introduction}
\bigskip
P.~L.~Ul'yanov\cite{Ul} raised the following question. Is it true that for any 
$2\pi$-periodic continuous function $f$ there is a uniformly convergent
rearrangement of its trigonometric Fourier series? The problem is still open. 

Let $\Bbb T=\Bbb R/2\pi\Bbb Z$, $C(\Bbb T)$ be the space of the continuous
functions $f:\Bbb T\to\Bbb C$, $\|f\|$ be the uniform norm of 
$f\in C(\Bbb T)$. We associate with every function $f\in C(\Bbb T)$ its 
Fourier series in a complex form
$$f\sim\sum_{k\in\Bbb Z}c_ke^{ikx}$$
and in a real form
$$f\sim\sum_{k=0}^\infty A_k(x),\quad A_k(x)=d_k\cos(kx+\phi_k).$$
Observe that $A_k(x)=c_ke^{ikx}+c_{-k}e^{-ikx}$.
It is easy to see that if Ul'yanov's conjecture is true for the series
in a real form (that is, there is a permutation $\sigma$ of $\Bbb N$
such that $\|f-d_0-\sum_{k=1}^n A_{\sigma(k)}\|\to0$ as $n\to\infty$) then it 
is also true for the series in a complex form because for $n\to\infty$
$$\left\|f-d_0-\sum_{k=1}^n \left(c_{\sigma(k)}e^{i\sigma(k)x}+c_{{-\sigma(k)}}
e^{-i\sigma(k)x}\right)\right\|\to0.$$

Sz.Gy.~~R\'ev\'esz\cite{R, R2} proved that for any $f\in C(\Bbb T)$ there is a 
rearrangement of its trigonometric Fourier series such that some subsequence of
the sequence of partial sums of a rearranged series converges to $f$
uniformly. Due to this result, Ul'yanov's conjecture is equivalent to the 
following: there is an absolute constant $C>0$ such that for any 
trigonometric polynomial (with a zero constant term) 
$\sum_{k=1}^n A_k(x)$ there is a permutation $\sigma:\{1,\dots,n\}\to
\{1,\dots,n\}$ such that for $m=1,\dots,n$
$$\left\|\sum_{k=1}^m A_{\sigma(k)}(x)\right\|
\le C\left\|\sum_{k=1}^n A_k(x)\right\|.$$ 
It is known that
$$\left\|\sum_{k=1}^m A_{k}(x)\right\|
\le C\log(n+1)\left\|\sum_{k=1}^n A_k(x)\right\|$$ 
(see \cite{Z, chapter~2, \S12}). Let 
$$\omega(f,\delta)=sup\Sb x,y\in\Bbb T\\|x-y|\le\delta\endSb
|f(x)-f(y)|$$
be the modulus of continuity of $f$. By Dini---Lipschitz theorem
\cite{Z, chapter~2, \S10},
if $\omega(f,\delta)=o(1/\log1/\delta)$ as $\delta\to0$ then
the Fourier series of $f$ converges to $f$ uniformly. Moreover, the condition
on $\omega(f,\delta)$ is sharp and cannot be replaced by
$\omega(f,\delta)=O(1/\log1/\delta)$\cite{Z, chapter~8, \S2}.

The author\cite{K, K2} proved the following results.
\proclaim{Theorem 1} For any trigonometric polynomial 
$\sum_{k=1}^n A_k(x)$ there is a permutation 
$\sigma:\{1,\dots,n\}$%\newline 
$\to\{1,\dots,n\}$ such that for $m=1,\dots,n$
$$\left\|\sum_{k=1}^m A_{\sigma(k)}(x)\right\|
\le C\log\log(n+2)\left\|\sum_{k=1}^n A_k(x)\right\|.$$ 
\endproclaim
\proclaim{Theorem 2} Let $f\in C(\Bbb T)$ and 
$\omega(f,\delta)=o(1/\log\log1/\delta)$ as $\delta\to0$. Then
there is a permutation $\Bbb N\to\Bbb N$ such that 
$$\left\|f-d_0-\sum_{k=1}^n A_{\sigma(k)}(x)\right\|
\to0\quad(n\to\infty).$$
\endproclaim
Theorem 2 follows from Theorem 1 by using Theorem 5 from \cite{R}. 

To approach Ul'yanov's conjecture, one can try to prove that there 
is an absolute constant $C>0$ such that for any 
trigonometric polynomial (with a zero constant term) 
$\sum_{k=1}^n A_k(x)$ and for any $m\le n$ there is an injection 
$\sigma:\{1,\dots,m\}\to\{1,\dots,n\}$ such that 
$$\left\|\sum_{k=1}^m A_{\sigma(k)}(x)\right\|
\le C\left\|\sum_{k=1}^n A_k(x)\right\|.$$ 
I cannot prove this either.

\proclaim{Theorem 3} For any trigonometric polynomial
$\sum_{k=1}^n A_k(x)$ and for any $m\le n$ there is a set  
$K\subset\{1,\dots,n\}$ such that $|K|=m$ and 
$$\left\|\sum_{k\in K} A_k(x)\right\|
\le C\log\log\log(n+20)\left\|\sum_{k=1}^n A_k(x)\right\|.$$ 
\endproclaim

\proclaim{Theorem 4} Let $f\in C(\Bbb T)$,
$$f\sim\sum_{k=0}^\infty A_k(x),\quad A_k(x)=d_k\cos(kx+\phi_k),$$
and $d_k=O(k^{-1/2})$. Then there is a permutation $\Bbb N\to\Bbb N$ such that 
$$\left\|f-d_0-\sum_{k=1}^n A_{\sigma(k)}(x)\right\|
\to0\quad(n\to\infty).$$
\endproclaim
In particular, Theorem 4 works if the sequence $\{|d_k|\}$
is nonincreasing. Note that, by theorem of Salem\cite{S}, 
there exists an even continuous function such that 
its Fourier series diverges at $x=0$
and the sequence $\{|d_k|\}$ is nonincreasing, 

By $C, C', C_1, C_2,\dots$ we denote positive constants
Let $[u]$ and $\{u\}$ be the integral and the fractional part 
of a real number $u$, respectively.
\bigskip
\centerline{\S2. Proof of Theorem 3}
\bigskip
Let $n\in\Bbb N$, $T$ be a trigonometric polynomial,
$$T(x)=\sum_{k=1}^n A_k(x)=\sum_{k=1}^n d_k\cos(kx+\phi_k).$$
We use the following lemmas from \cite{K2}.
\proclaim{Lemma 1} Let $\|T\|\le1$, $l\in\Bbb N$, $j\in\Bbb Z$,
$K_{l,j}=\{k:\ 1\le k\le n,\ k\equiv\pm j(\bmod l)\}$. Then 
$$\left\|\sum_{k\in K_{l,j}}A_k\right\|\le2.$$
\endproclaim
\proclaim{Lemma 2} Let $\|T\|\le1$. Then there exists an odd prime
$p\le2\log^3(n+3)$ such that
$$\sum\Sb k_1\neq k_2\\k_1\equiv k_2(\bmod p)\endSb |d_{k_1}|^2|d_{k_2}|^2
\le\frac{C_1}{\log^2(n+1)}.\tag1$$
\endproclaim
\proclaim{Lemma 3} Let $p$ be a prime satisfying (1), $j\in\Bbb Z$,
$K_{p,j}=\{k:\ 1\le k\le n,\ k\equiv\pm j(\bmod p)\}$, $N_j=|K_{p,j}|$. 
Then there exists a bijection $\tau:\{1,\dots,N_j\}\to K_{p,j}$
such that for any $m=1,\dots,N_j$ the inequality
$$\left\|\sum_{j=1}^m A_{\tau(j)}\right\|\le C_2(1+\|T\|)$$
holds.
\endproclaim

In the proof of Theorem 3  we assume that $n$ is sufficiently large and
and $\|T\|\le1$. We can also assume that $m\le n/2$; otherwise we can take 
the complement to a set constructed for $n-m<n/2$ instead of $m$. 
Also, it is sufficient to construct a set $K'\subset\{1,\dots,n\}$ such that 
$|K'|=m'$ for some $m'\le m$, $m-m'\le 0.2n/\log^3 n$, and
$$\left\|\sum_{k\in K'} A_k(x)\right\|\le C'\log\log\log n.$$ 
Indeed, take an odd prime $p\le 2\log^3 (n+3)$ satisfying Lemma 2.
Define the sets $K_{p,j}$ as in Lemma 3. Since $|K'|\le n/2$, we can find $j$ 
so that
$$|K_{p,j}\setminus K'|\ge(n-|K'|)/p\ge n/(4\log^3 (n+3))\ge 0.2n/\log^3 n$$
provided that $n\ge20$. Applying Lemma 3 to the polynomial
$$\sum_{k\in\{1,\dots,n\}\setminus K'}A_k,$$
we can define the set $K$ as $K'\cup\{\tau(1),\dots,\tau(m)\}$
where $m$ is such that 
$$\{\tau(1),\dots,\tau(m)\}\setminus K'=m-m'.$$
By the above arguments we can consider that 
$m>0.2n/\log^3 n$; otherwise, we take $m'=0$ and $K'=\emptyset$. 

We shall use the following known fact.
\proclaim{Lemma 4} For any real $\alpha\in(0,1]$ there exist positive integers 
$l_1,l_2,\dots,$ such that for any positive integer $s$
$$0<\alpha-\sum_{j=1}^s\frac1{l_j}\le2^{-2^{s-1}}.\tag2$$
\endproclaim
\demo{Proof of Lemma 4} We construct $l_s$ consequently:
$$l_s=\min\{l:\alpha-\sum_{j=1}^{s-1}\frac1{l_j}-\frac1l>0\}.$$
The inequalities (2) can be checked by induction on $s$. The proof of the
first inequality is straightforward. The induction base for the second 
inequality holds: $\alpha-1/l_1\le1/2$.

By induction supposition (2), we have 
$$l_{s+1}-1\ge2^{2^{s-1}}.$$
Also, by the definition of $l_{s+1}$,
$$\alpha-\sum_{j=1}^{s}\frac1{l_j}-\frac1{l_{s+1}-1}\le0.$$
Therefore,
$$\alpha-\sum_{j=1}^{s+1}\frac1{l_j}\le\frac1{l_{s+1}-1}-\frac1{l_{s+1}}
<\frac1{(l_{s+1}-1)^2}\le2^{-2^s},$$
and (2) is established for $s+1$. Lemma 4 is proved.
\enddemo

Take $s=[2\log\log\log n]$. Note that for sufficiently large $n$ we have 
$$2^{-2^{s-1}}\le 0.05/\log^3 n.\tag3$$ 
One can try to define the numbers $l_1,\dots,l_s$ by Lemma 3
with $\alpha$ close to $m/n$ and to take, for example,
$$K'=\bigcup_{j=1}^s K_j,\quad K_j=\{k\equiv\pm 1(\bmod 2l_j)\},$$
By Lemma 1,
$$\left\|\sum_{k\in K_j}A_k\right\|\le2$$
and $\sum_j|K_j|$ is close to $m$. However, the sets $K_j$
might have common points, and in general we cannot give good estimates
for $\left\|\sum_{k\in K'}A_k\right\|$ and for $|K'|$. We show how to
correct the construction.

Let $l_0=[5\log\log\log n]$, $\gamma=l_0m/n-0.1/\log^3 n$,
$g=[\gamma]$, $\alpha=\{\gamma\}$. Note that $g\ge0$  because of our 
supposition $m>0.2n/\log^3 n$. Take the numbers $l_1,\dots,l_s$ 
in accordance with Lemma 4 and define 
$$K'=\bigcup_{j=1}^\gamma K_j\cup\bigcup_{j=1}^s K_j',$$
where $K_j=\{k\equiv\pm j(\bmod 2l_0)\}$,
$K_j'=\{k\pm (g+j)\equiv0(\bmod 2ll_j)\}$.
Note that the residues classes $\pm j(\bmod 2l)$
$(j=1,\dots,\gamma+s)$, are all distinct since $\gamma+s\le l_0/2+s<l_0-1$. 
Therefore, the sets $K_j$, $K_j'$ are pairwise disjoint. Further, 
by Lemma 1,   
$$\left\|\sum_{k\in K_j}A_k\right\|\le2,\quad
\left\|\sum_{k\in K_j'}A_k\right\|\le2.$$
Hence,
$$\left\|\sum_{k\in K'}A_k\right\|\le2(g+s)\le10\log\log\log n.$$
Also, it is not difficult to check that
$$||K_j|-n/l_0|\le1, ||K_j'|-n/(l_0l_j)|\le1.$$
Therefore,
$$|K'|=n\gamma/l_0+\sum_{j=1}^s n/(l_0l_j)+O(\log\log\log n).$$

Taking into account (2) and (3), we get
$$n\gamma/l_0+\sum_{j=1}^s n/(l_0l_j)\le m-0.1n/\log^3 n,$$
$$n\gamma/l_0+\sum_{j=1}^s n/(l_0l_j)\ge m-0.1n/\log^3 n
-0.05n/\log^3 n.$$
Combining three last inequalities, we obtain
$$m\ge |K'|\ge m-0.2n/\log^3 n,$$
as required. This completes the proof of Theorem 3.
\bigskip
\centerline{\S3. Spencer's theorem and its corollaries}
\bigskip
Let $u$ be a vector: $u=(u^1,\dots,u^n)\in\Bbb R^n$. Denote
$|u|_\infty=\max_k|u^k|$.\newline  
J.~Spencer\cite{Sp} actually proved the following theorem.
\proclaim{Theorem A} Let $r\le n$ be positive integers, $u_j\in\Bbb R^n$, 
$|u_j|_\infty\le1$. Then for some choice of signs
$$|\pm u_1\pm\dots\pm u_r|_\infty\le C_3(r\log(2n/r))^{1/2}.$$ 
\endproclaim
\proclaim{Corollary 1}
Let $r\le n$ be positive integers and 
$K\subset\{1,\dots,n\}$, $|K|=r$. Consider a trigonometric polynomial
$$\sum_{k\in K}A_k(x),\quad A_k(x)=d_k\cos(kx+\phi_k).$$ 
Then there are sets $K_+\subset K$ and $K_-\subset K$ such that
$$K_+\cup K_-=K,\quad K_+\cap K_-=\emptyset,\quad |K_+|=|K|/2\tag4$$
and
$$\left\|\sum_{k\in K_+}A_k- \sum_{k\in K_-}A_k\right\|
\le C_4(r\log(2n/r))^{1/2}\max_{k\in K}|d_k|.\tag5$$ 
\endproclaim
\demo{Proof} Denote $d=\max_{k\in K}|d_k|$. We apply Theorem A to the vectors 
$u_k\in\Bbb R^{10n+1}$, $k\in K$, defined as
$$\gather
u_k=(\Re(A_x(\pi l/(5n))/d)_{l=0,\dots,2n-1},\\
\Im(A_x(\pi l/(5n))/d)_{l=0,\dots,2n-1},1).
\endgather$$
Then there exist numbers $\sigma_k=\pm1\ (k\in K)$
such that
$$\left\|\sum_{k\in K}\sigma_k A_k\right\|
\le 3C_3(r\log((20n+2)/r))^{1/2}d\tag6$$
and 
$$\left|\sum_{k\in K}\sigma_k\right|\le C_3(r\log((20n+2)/r))^{1/2}.\tag7$$
For the proof of (6) we use that for any trigonometric polynomial $T$ of 
order $n$
$$\|T\|\le3\max_{l=0,\dots,10n-1}|T(\pi l/(5n))|$$
(see, for example, \cite{Kl}).
Without loss of generality we can assume that $\sum_{k\in K}\sigma_k\le0$.
Take $K_+'=\{k\in K:\ \sigma_k=1\}$, $K_-'=\{k\in K:\ \sigma_k=-1\}$.
We have
$$2|K_+'|=|K|+\sum_{k\in K}\sigma_k\le2[|K|/2].$$
Take an arbitrary set $K_1\subset K_-$ such that $|K_-|=[|K|/2]-|K_+'|$.
By (7), $|K_1|\le C_3(r\log((20n+2)/r))^{1/2}/2$. Hence,
$$\left\|\sum_{k\in K_1}A_k\right\|
\le C_3(r\log((20n+2)/r))^{1/2}d.\tag8$$
Denote $K_+=K_+'\cup K_1$, $K_-=K_-'\setminus K_1$. The conditions (4)
are satisfied. By (6) and (8) we get
$$\left\|\sum_{k\in K_+}A_k- \sum_{k\in K_-}A_k\right\|
\le 3.5C_3(r\log((20n+2)/r))^{1/2}d.$$ 
Therefore, (5) also holds, and Corollary 1 is proved.
\enddemo

\proclaim{Corollary 2}
Let $r\le n$ be positive integers and 
$K\subset\{1,\dots,n\}$, $|K|=r$. Consider a trigonometric polynomial
$$\sum_{k\in K}\alpha_kA_k(x),\quad A_k(x)=d_k\cos(kx+\phi_k),$$
where $\alpha_k$ are real numbers. Then there are numbers 
$\beta_k\in\{[\alpha_k],[\alpha_k]+1\}$ such that 
$$\left\|\sum_{k\in K}\alpha_kA_k-\sum_{k\in K}\beta_kA_k\right\|
\le C_4(r\log(2n/r))^{1/2}\max_{k\in K}|d_k|.$$ 
\endproclaim
In fact, the deduction of Corollary 2 from Corollary 1 is 
exhibited in \cite{Kl}.

\proclaim{Corollary 3} Let $r, n$ be positive integers, $r\le n/5$ and 
$K\subset\{1,\dots,n\}$, $|K|=r$. Consider a trigonometric polynomial
$$\sum_{k\in K}A_k(x),\quad A_k(x)=d_k\cos(kx+\phi_k).$$ 
Then there exists  a bijection $\sigma:\{1,\dots,r\}\to K$
such that for any $m=1,\dots,r$ the inequality
$$\gather
\left\|\sum_{j=1}^m A_{\sigma(j)}-\frac mr \sum_{k\in K}A_k\right\|\\
\le (4C_4+4)(r\log(2n/r))^{1/2}\max_{k\in K}|d_k|\tag9
\endgather$$
holds. 
\endproclaim
\demo{Proof} Denote $d=\max_{k\in K}|d_k|$. We fix $n$ and use induction on 
$r$. If $r\le8$ then we take an arbitrary bijection $\sigma$. For any
$m\le r$ we have
$$\gather
\left\|\sum_{j=1}^m A_{\sigma(j)}-\frac mr \sum_{k\in K}A_k\right\|
\le md+\frac mr(rd)\le2md\\
\le2rd= (2r)^{1/2}(2r)^{1/2}d\le4(r\log(2n/r))^{1/2}d,
\endgather$$
and (9) holds. Let us assume that $9\le r\le n/5$ and that the statement of 
the corollary is satisfied for all $r'<r$.

By Corollary 1, we split the sets $K$ into the sets $K_+$ and $K_-$.  
The inequality (5) can be rewritten as
$$\left\|\sum_{k\in K_+}A_k-\frac12 \sum_{k\in K}A_k\right\|
\le\frac{C_4}2(r\log(2n/r))^{1/2}d.$$ 
We have
$$\gather
\left\|\sum_{k\in K_+}A_k-\frac{[r/2]}r \sum_{k\in K}A_k\right\|
\le\frac{C_4}2(r\log(2n/r))^{1/2}d\\ 
+\left(\frac12-\frac{[r/2]}r\right)\left\|\sum_{k\in K}A_k\right\|
\le\frac{C_4}2(r\log(2n/r))^{1/2}d\\ 
+\frac1{2r}(rd)=\frac{C_4}2(r\log(2n/r))^{1/2}d+d/2\\
\le\frac{C_4+1}2(r\log(2n/r))^{1/2}d.\tag10
\endgather$$

By the induction supposition, there exist bijections 
$\sigma_+:\{1,\dots,[r/2]\}\to K_+$ and
$\sigma_+:\{1,\dots,r-[r/2]\}\to K_-$ such that for any $m\le[r/2]$
$$\gather
\left\|\sum_{j=1}^m A_{\sigma_+(j)}-\frac m{r_1} \sum_{k\in K_+}A_k\right\|\\
\le (4C_4+4)(r_1\log(2n/r_1))^{1/2}d,\quad r_1=[r/2],\tag11
\endgather$$
and for any $m\le r-[r/2]$
$$\gather
\left\|\sum_{j=1}^m A_{\sigma_-(j)}-\frac m{r_1} \sum_{k\in K_-}A_k\right\|\\
\le (4C_4+4)(r_1\log(2n/r_1))^{1/2}d,\quad r_1=r-[r/2].\tag12
\endgather$$

We take $\sigma(j)=\sigma_+(j)$ for $j\le[r/2]$ and
$\sigma(j)=\sigma_-(r+1-j)$ for $j>[r/2]$. If $m\le[r/2]$
then we have, by (10) and (11),
$$\gather
\left\|\sum_{j=1}^m A_{\sigma(j)}-\frac mr \sum_{k\in K}A_k\right\|
\le\left\|\sum_{j=1}^m A_{\sigma_+(j)}-\frac m{r_1} 
\sum_{k\in K_+}A_k\right\|\\
+\left\|\frac m{r_1} \sum_{k\in K_+}A_k-\frac mr \sum_{k\in K}A_k\right\|\\
\le\left\|\sum_{j=1}^m A_{\sigma_+(j)}-\frac m{r_1} 
\sum_{k\in K_+}A_k\right\|\\
+\left\|\sum_{k\in K_+}A_k-\frac{[r/2]}r \sum_{k\in K}A_k\right\|\\
\le (4C_4+4)(r_1\log(2n/r_1))^{1/2}d\\
+\frac{C_4+1}2(r\log(2n/r))^{1/2}d,\quad r_1=[r/2].\tag13
\endgather$$
Further, for $r_1=[r/2]$ we have
$$\gather
(r_1\log(2n/r_1))^{1/2}\le\left(\frac r2\log(2n/r\times9/4)\right)^{1/2}\\
<\left(\frac r2\times\frac32\log(2n/r)\right)^{1/2}
<\left(\frac34r\log(2n/r)\right)^{1/2}\\
<\frac78(r\log(2n/r))^{1/2}.
\endgather$$
Substituting the last inequality into (13) we get the required
$$\left\|\sum_{j=1}^m A_{\sigma(j)}-\frac mr \sum_{k\in K}A_k\right\|
\le(4C_4+4)(r\log(2n/r))^{1/2}d.$$
If $m>[r/2]$, then, similarly to (13), e have 
$$\gather
\left\|\sum_{j=1}^m A_{\sigma(j)}-\frac mr \sum_{k\in K}A_k\right\|\\
=\left\|\sum_{j=1}^{r-m} A_{\sigma_-(j)}-\frac {r-m}r 
\sum_{k\in K}A_k\right\|\\
\le (4C_4+4)(r_1\log(2n/r_1))^{1/2}d\\
+\frac{C_4+1}2(r\log(2n/r))^{1/2}d,\quad r_1=r-[r/2].\tag14
\endgather$$
For $r_1=[r/2]$ we have
$$\gather
(r_1\log(2n/r_1))^{1/2}\le\left(\frac{5r}9\log(2n/r\times2)\right)^{1/2}\\
<\left(\frac{5r}9\times\frac43\log(2n/r)\right)^{1/2}
<\left(\frac34r\log(2n/r)\right)^{1/2}\\
<\frac78(r\log(2n/r))^{1/2}.
\endgather$$
and after substitution of the last inequality into (14) we complete the proof 
of Corollary 3.
\enddemo
\bigskip
\centerline{\S4. Proof of Theorem 4}
\bigskip
We use Vall\'ee Poussin sums defined for positive integers $n>m$ as
$$V_{m,n}(x)=\sum_{k=0}^m A_k(x)+\sum_{k=m+1}^n\frac{n-k}{n-m}A_k(x).$$
It is known that for any $f\in C(\Bbb T)$ there is a function 
$n:\Bbb N\to\Bbb N$ such that $n(m)>m$ for all $m$, 
$\lim_{m\to\infty}n(m)/m=1$ and $\lim_{m\to\infty}\|V_{m,n}-f\|=0$
(this follows, for example, from \cite{D} or from \cite{St}).  
We define the increasing sequence of positive integers 
$\{N_\lambda\}_{\lambda\in\Bbb N}$ by $N_1=1$, 
$N_{\lambda+1}=n(N_{\lambda})$ for $\lambda\ge1$. 

We fix $\lambda\ge1$, take $m=N_{\lambda}$, $n=N_{\lambda+1}$ 
and use Corollary 2 for $K_\lambda=\{m+1,\dots,n\}$, 
$\alpha_k=\frac{n-k}{n-m}$. We find that there are numbers 
$\beta_k\in\{0,1\}$, $k\in K$, such that
$$\gather
\left\|V_{m,n}-\sum_{k=0}^m A_k-
\sum_{k\in K}\beta_kA_k\right\|\\
\ll(((n-m)/n)\log((2n)/(n-m)))^{1/2}\to0\quad(\lambda\to\infty).\tag15
\endgather$$
Also, by the choice of the sequence $\{N_\lambda\}\}$, we have\newline
$\lim_{\lambda\to\infty}\|V_{m,n}-f\|=0$. Therefore, 
denoting 
$$L_{\lambda}=\{1,\dots,m\}\cup\{k\in K_\lambda:\ \beta_k=1\}$$
we get
$$\left\|f-d_0-\sum_{k\in L_{\lambda}}A_k\right\|\to0\quad(\lambda\to\infty).
\tag16$$

To complete the proof, it is enough, by (16), to find a good permutation
of the terms of the polynomials
$$\sum_{k\in L_{\lambda+1}\setminus L_\lambda}A_k.$$
We construct a permutation in such a way that the numbers from 
$L_\lambda\setminus L_{\lambda-1}$ precede the numbers from  
$L_{\lambda+1}\setminus  L_\lambda$ for all $\lambda$ for all 
$\lambda\in\Bbb N$; we consider that $L_0=\emptyset$.
The permutation can be constructed by Corollary 3, the partial sums
can be estimated similarly to (16), and we are done. 
\bigskip
\centerline{REFERENCES}
\bigskip

[D] V.~Damen, Best approximations and de la Vall\'ee\newline 
Poussin sums (Russian), Mat. Zametki, 23 (1978), pp.~671--683. 

[K] S.~V.~Konyagin, On rearrangements of trigonometric Fourier series,
(Russian), Vsesoyuznaya shkola ``Teoriya priblizheniya funkciy''.
Tezisy dokladov, p.~80. Kiev, 1989.

[K2] S.~V.~Konyagin, On uniformly converging rearrangements of trigonometric 
Fourier series, (Russian), Metric theory of functions and related problems 
in analysis (Russian), pp.~101--111, Izd. Nauchno-Issled. Aktuarno-Finans. 
Tsentra (AFTs), Moscow, 1999. 

[Kl] M.~N.~Kolountzakis, On nonnegative cosine polynomials with nonnegative 
integral coefficients, Proc.~AMS, 120 (1994), pp.~157--163. 

[R] Sz.Gy.~R\'ev\'esz, Rearrangement of Fourier series,\newline
J.~Appr.~Theory, 60 (1990), pp.~101--121

[R2] Sz.Gy.~R\'ev\'esz, On the convergence of Fourier series 
of U.A.P. functions, J.~Math.~Anal.~Appl., 151 (1990),
pp.~308--317.
 
[S] R.~Salem, On a problem of Smithies, Indag.~Math., 16 (1954),
pp.~403--407.

[Sp] J.~Spencer, Six standard deviation suffice, Trans. Amer.~Math.~Soc.,
289 (1985), pp.~679--706.

[St] S.~B.~Stechkin, On the approximation of periodic functions by de la 
Vall\'ee Poussin sums, Anal.~Math., 4 (1978), pp.~61--74.

[Ul] P.~L.~Ul'yanov, Solved and unsolved problems in the theory of 
trigonometric and orthogonal series, (Russian), Uspehi Mat. Nauk,
19 (1964), N.~1(115), pp.~3--69. 

[Z] A.~Zygmund, Trigonometric series, v.~1. Cambridge, The University Press,
1959.
\enddocument
\end